\documentclass[letterpaper, 10 pt, conference]{ieeeconf}  
\usepackage{graphicx}

\IEEEoverridecommandlockouts                              

\title{\LARGE \bf
Predicting and controlling the dynamics of infectious diseases
}

\author{Robin J. Evans$^{1}$ and Musa Mammadov$^{1}$
\thanks{$^{1}$University of Melbourne, Parkville VIC 3010, Australia
        {\tt\small robinje@unimelb.edu.au, mmammadov@unimelb.edu.au}}%
}

\begin{document}

\maketitle
\thispagestyle{empty}
\pagestyle{empty}

\begin{abstract}

This paper introduces a new optimal control model to describe and control the dynamics of infectious diseases.
In the present model, the average time of isolation (i.e. hospitalization) of infectious population is the main time-dependent parameter that defines the spread of infection.
All the preventive measures aim to decrease the average time of isolation under given constraints.

\end{abstract}

\section{INTRODUCTION}

The outbreak of Ebola Virus Disease (EVD) in West Africa 2014
revealed many challenges in predicting and controlling the spread of infectious diseases.
These challenges are partly related to the mathematical modeling of the dynamics of the epidemic.
Providing accurate predictions appeared to be extremely difficult.

To address these challenges, several new models have been suggested each providing quite different results, for example  \cite{althaus2014estimating,browne2014model,chowell2004basic,chowell2014transmission,nishiura2014early,rivers2014modeling,periods2014ebola}.
We also note that different aspects of possible control have been intensively studied in the literature including distribution strategies for vaccination and antibiotic programs \cite{Andrews2011} as well as travel restrictions \cite{Epstain2007}.

In this paper we concentrate on the development of  models that are well suited to the control of an outbreak.
The most commonly studied models in this area deal with temporal networks \cite{Li2006,Masuda2013,Salate2010} where several different models have been suggested.

In \cite{Lekone2006} the 1995 Ebola outbreak in Congo is considered using an SEIR model whereby control intervention, performed at time $t_*,$ is described by the transmission coefficient $\beta(t)$ defined by $\beta (t) = \beta$ if $t < t_*$ and $\beta (t) = \beta \exp(-q (t-t_*))$ for $t \ge t_*,$ where $\beta$ is the initial transmission rate that would remain stable without the intervention.

Our paper models and studies an alternative  control
mechanism to decrease in the transmission of infection based on the model developed in \cite{F1000} where
the transmission of infection depends mainly on two parameters: the transmission rate $\beta$
and the average time of isolation $\tau.$  In contrast to \cite{Lekone2006}
we assume that the transmission rate $\beta$ does not change over the whole period under consideration.
The main dynamic parameter in our model is $\tau$ and all the intervention measures
are directed at decreasing $\tau$ and consequently  reducing the spread of infection.

The average time to hospitalization can be used for the average time of isolation $\tau,$  although the isolation of infectious population is not exactly the same as ``hospitalization".
Clearly, as the disease progresses hospitals become short on beds (as well as staff and supplies)
required to isolate and treat all newly infected individuals.
As a result, the number of infected population can grow exponentially.
This was the case for the Ebola virus epidemic in West Africa (Guinea, Sierra Leone and Liberia) where the spread of infection
was highly dangerous during June-November 2014 when the capacity for treating Ebola patients was insufficient.
It was reported \cite{Cooper} that during this period
``... many clinics and hospitals in all three of the countries worst hit by Ebola have effectively been shut down".

After rapidly building new infrastructure and increasing the capacity of beds
the outbreak slowed down significantly.
Starting from January 2015, the epidemic has moved to the ending phase that involves
ensuring ``capacity for case finding, case management, safe burials and community engagement" (\cite{WHO_data} - WHO, Ebola Situation Report, 28 Jan 2015).
Note that in \cite{Drake} the hospitalization rate was the parameter showing the greatest change.

Addressing these issues, this paper suggests new mathematical models that can be used to
increase the efficiency of available resources.
The main goal here is to keep the model as simple as possible and, at the same time, to have measurable control variables.
Note that there are many useful control measures that have been intensively studied by introducing more ``detailed" mathematical models,
however such models have less predictive capabilities (due to overfitting). Prediction is crucial when considering future planning periods.

The main component in the suggested model is the optimal distribution of bed capabilities
across countries/regions. This is a very important and difficult problem
that requires an accurate prediction of the dynamics of infected population in each region.
For example, evaluating the situation of the Ebola outbreak, WHO's Ebola Situation Report on 14 Jan 2015 \cite{WHO_data}
notes that ``Each of the intense-transmission countries has sufficient capacity to isolate
and treat patients, with more than 2 treatment beds per reported
confirmed and probable case. However, the uneven geographical distribution
of beds and cases, and the under-reporting of cases,
means that not all EVD cases are isolated in several areas."

\section{MODEL}

In \cite{F1000} a new model is introduced to study the dynamics of epidemics by considering
the average time for isolation (denoted by $\tau$) of infectious population as a time-dependent parameter.
This model is derived from the well studied $SIR$ (Susceptible-Infectious-Recovery) model
(e.g. \cite{li2014sir}) and is similar to models based on transmission rates from infectious population
at different generations (e.g. \cite{nishiura2014early}).

The use of time-dependent parameter $\tau$ enables the analysis of future scenarios by considering
possible changes in $\tau.$ In this paper we extend this approach
by developing practical and efficient optimal control models.

We denote by $x(t)$ the number of infected cases at $t \in \{1,2, \cdots, T\}$ (in days).
Assuming that the natural death rate of population ($\mu$) is zero, the equation for $x(t)$ is as follows (see \cite{F1000} for more details)
\begin{equation}\label{x}
     x(t+1) = \beta  \, \sum_{i = 0}^{\tau-1} (1-\alpha \omega(i)) \, x(t-d - i).
\end{equation}
Here $\alpha$ is the death rate due to disease;
$d$ is the average latent period (in days) for infected individuals to become infectious;
$\tau$ is the average infectiousness period (in days); it is the average time required for isolation (time to hospitalization); and
$\beta$ is the transmission rate.
Moreover, $\omega$ is a gamma (cumulative) distribution function (with p.d.f - $\omega_p$) for deaths due to disease.
The fraction $(1-\alpha \omega(i))$ in this case represents
the proportion of remaining infected cases $x(t-d - i)$ after $d+i$ days.

The sum $$I_a(t) = \sum_{i = 0}^{\tau-1} (1-\alpha \omega(i)) \, x(t-d - i)$$
defines the number of "active" infectious population at time $t;$
it represents the number of infectious population that are not yet isolated and therefore it is the only source of secondary infections.
(for the sake of simplicity we do not consider infections in hospitals and death ceremonies).
Then by setting $x(t) = \beta I_a(t)$ we obtain model (1) in  \cite{F1000} where $\mu$ (the natural death rate) in our case is zero.

\bigskip

The basic reproduction number $R$ is calculated by considering the stationary states in (\ref{x}):
\begin{equation}\label{R}
    R = \beta  \, [\tau - \alpha \sum_{i = 0}^{\tau-1}  \omega(i)].
\end{equation}

There are three main parameters in (\ref{x}) - $\alpha, \beta$ and $\tau.$
The results obtained in \cite{F1000} show that this model provides quite good approximation to the
total number infected cases and deaths during the current Ebola epidemic
if $\tau$ is a piecewise
constant function (in fact, constant over consequent subintervals with durations 2-3 months)
and the parameters $\alpha$ and $\beta$ are constant over the whole period.

These results help us to predict the dynamics of an infected population at future time intervals by keeping the values of $\alpha$ and $\beta$
unchanged (estimated from the previous period) and considering different possible changes in $\tau.$
In this case the major strategy of preventive intervention is the achievement of some decrease in
$\tau$ that according to (\ref{R}) is equivalent to decreasing the effective reproduction number.

This approach is implemented below by introducing an optimal control models where
the average time to hospitalization $\tau$ is the key variable.
According to the results of data fitting mentioned above, it is sufficient to let
$\tau$ be constant on quite long time intervals (months).\\

{\bf Control $\mathbf{\tau}$. ~ }
Therefore, we define $\tau(t)$ as a control variable by assuming that it is piece-wise constant with integer values (days). For the sake of simplicity let
 $$\mathbf{\tau}(t) \equiv \tau_{i} \in U, ~ \forall t \in (T_{j}, T_{j+1}], ~ j = 1, 2, \cdots, p.$$
It is reasonable to assume that
$U \doteq \{\tau_{min},\tau_{min}+1 ,\cdots, \tau_{max}\};$ where $\tau_{min} \ge 1$ is the minimal number of days required to isolate infectious population.\\

{\bf Trajectory $\mathbf{x}$. ~ }
Given control $\mathbf{\tau}$ we define trajectory $\mathbf{x} = x (t)$ as follows

\begin{equation}\label{xrt}
     x(t+1) = \beta  \, \sum_{i = 0}^{\tau (t)-1} (1-\alpha \omega(i)) \, x(t-d - i).
\end{equation}
In this formula the sum $\sum_{i = 0}^{\tau (t)-1} (1-\alpha \omega(i)) \, x(t-d - i)$ represents the
number of infectious individuals that are not yet isolated.
Considering the average length of hospital stay (in days), the number of hospitalized cases at $t$
can be calculated as
\begin{equation}\label{H}
  h(t: \mathbf{\tau},\mathbf{x}) = \sum_{i = \tau (t)}^{\sigma} (1-\alpha \omega(i)) \, x(t-d - i)
\end{equation}
Note that, recent studies (see for example \cite{Drake}) show that $(\sigma - \tau (t))$  is around 6.5 days.

\section{Data fitting}\label{fit}

In this section we provide some numerical experiments based on data from Guinea, Sierra-Leon and liberia.
We consider the cumulative number of infectious cases and deaths  denoted by
$C(t)$ and $D(t),$ respectively. They can be calculated as

\begin{equation}\label{C}
     C(t+1) =  \sum_{s = 0}^{t} x(s-d);
\end{equation}

\begin{equation}\label{D}
     D(t+1) =  \sum_{s = 0}^{t} \, \sum_{i = 0}^{n} \alpha \omega_p(i) \, x(s-d - i).
\end{equation}
Here  $\alpha \omega_p(i)$ is the death rate of infectious population in generation $x(s-d - i)$ and $n$ is a large number.
Parameters of the gamma distribution function $\omega(i)$ are taken from \cite{periods2014ebola} where
\begin{equation}\label{gamma}
    \omega_p(x) = \frac{b^a}{\Gamma(a)} x^{a-1} e^{-bx}, ~~ a = 10, ~ b = 1.3333
\end{equation}
with mean value 7.5. Moreover, we set $d=7$ and $n=35$ as in \cite{F1000}.

In the considered model (\ref{xrt}) there are only three
parameters $\alpha$ and $\beta$ (constants) and a piece-wise constant control function $\tau(t)$ that need to be
optimized to fit data - the total number of infectious population and deaths.
The aim here is to show that there exists a control $\tau(t)$ such that
the corresponding trajectory $x(t)$ fits data well.\\

\begin{figure}[thpb]
      \centering
\begin{tabular}{c}
\includegraphics[scale=0.3]{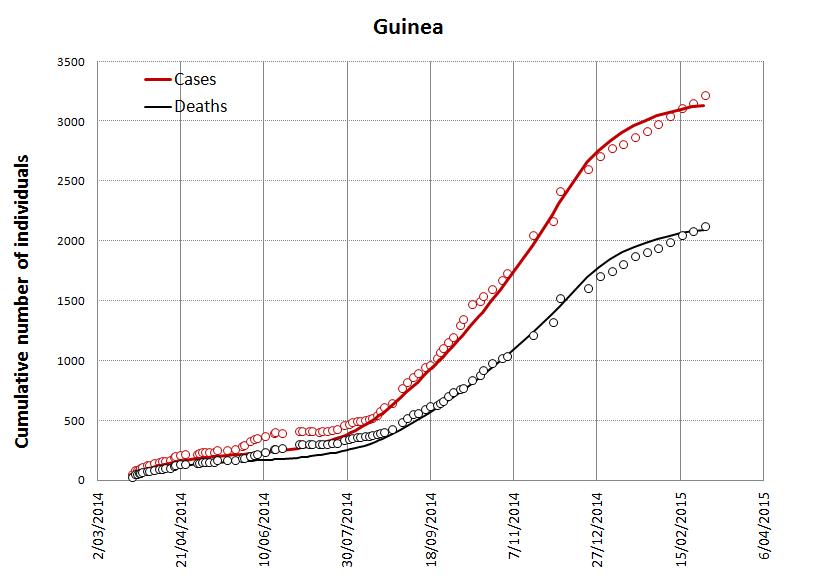} \\ \includegraphics[scale=0.3]{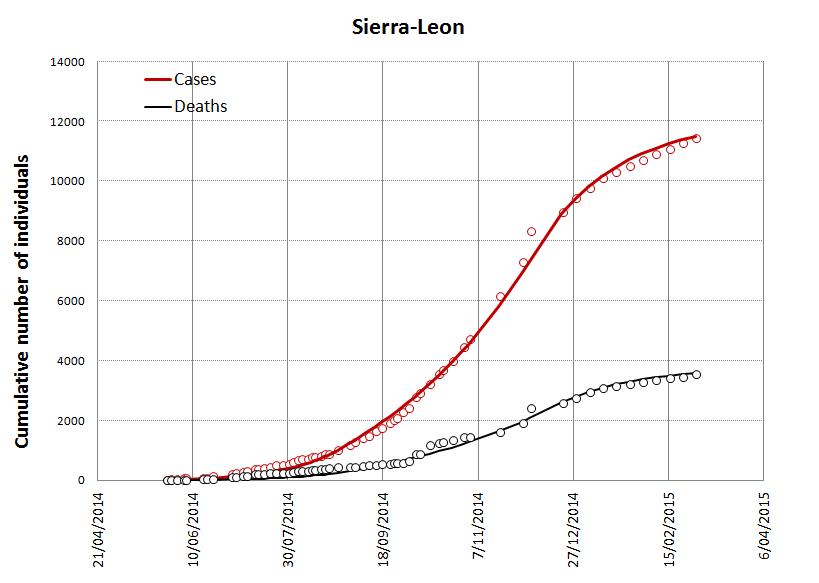} \\
\includegraphics[scale=0.3]{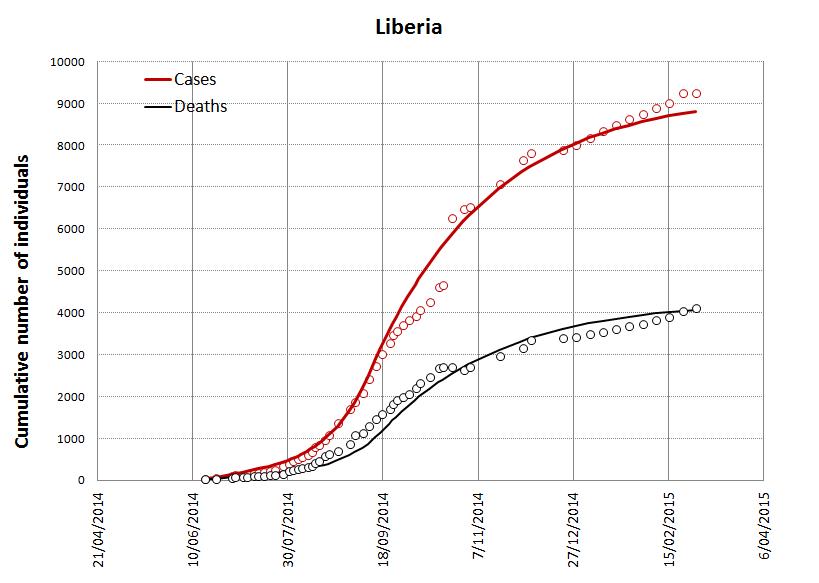} \\
\end{tabular}
\caption{The best fits for the cumulative numbers of infected cases and deaths in Guinea, Sierra Leone and
Liberia by considering three parameters $\alpha$, $\beta$ and 4 subintervals
with different values $\tau_k,$ $k=1,2,3,4$ (for the values see Table \ref{Table_R_k}).
The lines represent the best fits, red and black circles represent the data}
      \label{Fig02}
   \end{figure}

We consider three consequent intervals $\Delta_k = [T_k, T_{k+1}]$ $(k = 1,\cdots,4)$ for each country and find optimal values
$\alpha,$ $\beta$ and $\tau_k$ $(k = 1,\cdots,4)$ where $\tau(t) = \tau_k, ~ \forall t \in \Delta_k.$ The results are presented in Table \ref{Table_R_k}.
The last time point $T_5$ is 01-Mar-2015. The values of $T_1, T_2,T_3, T_4$ are as follows:
22-March, 23-May, 20-July and 04-Dec-2014 for Guinea;
27-May, 20-June, 20-August and 04-Dec-2014 for Sierra Leone; and
16-June, 20-July,  07-Sept and 04-Dec-2014 for Liberia.
Each interval $\Delta_k$ has its own reproduction number $R_k$ that defines the shape of the best fits
presented in Figure \ref{Fig02}.

Data were retrieved from the WHO website \cite{WHO_data} for the cumulative numbers of clinical cases (confirmed, probable and suspected) collected till 1 March 2015. The global optimization algorithm DSO in Global And Non-Smooth Optimization (GANSO) library \cite{Ganso, Mam-2005-Chapter} is applied for solving optimization problems in this section as well as in Section \ref{example}.

\begin{table}
\caption{{\small Results of best fits: the (effective) reproduction numbers $R_k$ and average times to hospitalization $\tau_k$ (in days) for different intervals $\Delta_k,$ $k = 1,2,3,4.$ The optimal values for $\alpha$ and $\beta$ are also provided; they are constant for a whole period}}\label{Table_R_k}
\medskip
\centering
\begin{tabular}{|c|cc|cccc|}
  \hline
 Country & $\alpha$ & $\beta $ & $R_1$ ($\tau_1$)  & $R_2$ ($\tau_2$)& $R_3$ ($\tau_3$)& $R_4$ ($\tau_4$) \\
\hline
 Gui.     & 0.66 & 0.265 & 0.79 (3)  & 1.31 (5) & 1.06 (4) & 0.79 (3)\\
 S.-L.    & 0.32 & 0.274 & 1.36 (5)  & 1.36 (5) & 1.09 (4) & 0.82 (3)\\
 Lib.     & 0.46 & 0.294 & 1.17 (4)  & 1.46 (5) & 0.88 (3) & 0.88 (3)\\
 \hline
\end{tabular}
\end{table}

The results obtained show that the estimated values of $\alpha$ and $\beta$ can be used for future time intervals while considering
$\tau$ as a dynamic parameter that defines the spread of infection.
This naturally leads to  optimal control problems that are considered in the next section.

\section{OPTIMAL DISTRIBUTION OF BED CAPACITIES}

Denote by $B(t)$ the number of beds (capacity of hospitals) at time $t.$ It is an increasing
and piece-wise constant function where jumps are at the points $T^b_{i}, ~ i = 1, 2, \cdots, q.$
We refer to \cite{Drake} (Fig 1) for an example of $B(t)$ in Liberia between June-September 2014.
Note that $q \ge p$; that is the number of points $T^b_{i}$ is larger than $T_{j}.$\\

We assume that the initial period of epidemic is long enough to estimate parameters $\alpha,$ $\beta$
and to provide some future scenarios depending on $\tau(t).$ \\

\subsection{Objective function}

Given process $(\mathbf{\tau},\mathbf{x})$ we define hospitalization rate as
\begin{equation}\label{HR}
  HospRate(t) = \frac{h(t: \mathbf{\tau},\mathbf{x})}{B(t)}.
\end{equation}

\bigskip
Note that the number of infectious cases (especially during high growth) may exceed the number of available beds, while it becomes quite low when the infection is slowing down. Taking into account this observation, we use
the hospitalization rate in the definition of the objective function (to be minimized) given by
$$
  x(T) + K \cdot \sum_{t=1}^{T}\frac{h(t: \mathbf{\tau},\mathbf{x})}{B(t)};
$$
where $K$ is constant.

\subsection{Constraints in terms of costs}

We consider the following three functions that will be used to formulate cost constraints.

\begin{itemize}
  \item $c_{B}(\Delta b)$ - costs associated to building $\Delta b$ additional beds;
  \item $c_{S}(h)$  - costs required for servicing $h$ patients in hospitals;
  \item $c_{I}(\Delta h)$ - costs required for $\Delta h$ infectious cases before hospitalization.
\end{itemize}

\bigskip
Assuming that $F(t)$ is the total available funds, we can formulate cost constraints as\\
$$
  c_{B}(\Delta b(t)) + c_{S}(h(t)) + c_{I}(\Delta h(t)) \le F(t), ~ \forall t.
$$

\bigskip
Here $\Delta b(t) = b(t+1)-b(t)$ and $\Delta h(t) = h(t+1: \mathbf{\tau},\mathbf{x}) - h(t: \mathbf{\tau},\mathbf{x}).$

\subsection{Multi regions}\label{multi}

Now consider $m$ regions/countries and assume there is an inter-transmission of infections between them.
Denote by $x_r(t)$ the number of infected cases in country $r$ and assume that the transmission of infection from
country $i$ to $r$ is given by the coefficient $\beta_{ir}$ (\cite{nishiura2014early}).
In this case, by considering corresponding average times of isolation $\tau_r(t)$
in each country $r,$ we have the following system
$$
  \begin{array}{cc}
    x_r(t+1) & = \beta_{1r}  \, \sum_{i = 0}^{\tau_1 (t)-1} (1-\alpha_1 \omega(i)) \, x_1(t-d - i) \\
             & + \beta_{2r}  \, \sum_{i = 0}^{\tau_2 (t)-1} (1-\alpha_2 \omega(i)) \, x_2(t-d - i) \\
             &  \cdots \\
             & + \beta_{mr}  \, \sum_{i = 0}^{\tau_m (t)-1} (1-\alpha_m \omega(i)) \, x_m(t-d - i) \\
  \end{array}
$$
where $\beta_{ir}$ is the infections generated from country $i.$ One would expect
$\beta_{rr}$ to be much larger than $\beta_{ir}$ for $i \neq r$ and
in a special case $\beta_{ir} = 0.$ The number of hospitalized population at time $t$ is defined by (\ref{H}); that is

$$
  h_r(t: \mathbf{\tau_r},\mathbf{x_r}) = \sum_{i = \tau_r (t)}^{\sigma} (1-\alpha_r \omega(i)) \, x(t-d - i), ~~ r= 1, \cdots, m.
$$

We assume that the death rates $\alpha_r$ and the coefficients $\beta_{ir}$ are
estimated from the initial data and they are constant over the whole period of the epidemic.
Then, we can consider the problem of optimal distribution of available new beds between regions formulated below.\\

{\bf Problem 1 (Optimal distribution of bed capacities):} Given initial data ($\alpha_r,$  $\beta_{ir}$ and $x_r(t), t \le 1$), total bed capacity $B(t)$ and future scenario for control functions $\tau_r(t), r = 1, \cdots, m,$ find increasing piece-wise constant functions $b_r(t), r = 1, \cdots, m,$ for the problem

$$
  {\rm Minimize:~}  \sum_{r=1}^{m} \left[ x_r(T) + K \cdot \sum_{t=1}^{T}\frac{h_r(t: \mathbf{\tau_r},\mathbf{x_r})}{b_r(t)} \right];
$$
$$
{\rm ~subject~to:~}  b_1(t) + \cdots + b_m(t)~ \le~ B(t), ~ \forall t.
$$

\bigskip
In the next problem we take into account the cost constraints:\\

{\bf Problem 2 (Optimal distribution of bed capacities under cost constraints):} Given initial data ($\alpha_r,$  $\beta_{ir}$ and $x_r(t), t \le 1$), total budget function $F(t)$ and future scenario for control functions $\tau_r(t), r = 1, \cdots, m,$ find increasing piece-wise constant functions $b_r(t), r = 1, \cdots, m,$ for the problem

$$
  {\rm Minimize:~}  \sum_{r=1}^{m} \left[ x_r(T) + K \cdot \sum_{t=1}^{T}\frac{h_r(t: \mathbf{\tau_r},\mathbf{x_r})}{b_r(t)} \right];
$$
$$
{\rm ~subject~to:~} b_1(t) + \cdots + b_m(t)~ \le~ B(t), ~ \forall t;
$$
$$
 c_{B}(\Delta B(t)) +  \sum_{r=1}^{m} \left[ c_{S}(h_r(t)) + c_{I}(\Delta h_r(t))\right] \le F(t), ~ \forall t.
$$

{\bf Feasible processes.~}
Given initial data ($\alpha_r,$  $\beta_{ir}$ and $x_r(t), t \le 1$) and $B(t)$) consider the trajectory
$\mathbf{x} = (\mathbf{x_1}, \cdots, \mathbf{x_m})$ corresponding to $\mathbf{\tau} = (\mathbf{\tau_1}, \cdots, \mathbf{\tau_m}).$ Denote also  $\mathbf{b} = (\mathbf{b_1}, \cdots, \mathbf{b_m}),$ where each $\mathbf{b_r}$
stands for the bed capacity function $b_r(t)$ in region $r.$

We call process $(\mathbf{\tau}, \mathbf{b}, \mathbf{x})$ feasible if all the constraints
of the problem under consideration hold and the hospitalization rates are less than 1;
that is, $$h_r(t: \mathbf{\tau_r},\mathbf{x_r}) \le b_r(t), ~ \forall t, ~ r = 1, \cdots, m.$$

The meaning of feasible processes can be explained as follows. If $(\mathbf{\tau}, \mathbf{b}, \mathbf{x})$ is feasible
then the number of required beds and the resources needed for
isolation are sufficient at every time point $t$ in order to keep the average times of isolation at level $\mathbf{\tau} = (\mathbf{\tau_1}, \cdots, \mathbf{\tau_m}).$
Thus, the corresponding effective reproduction numbers $\mathbf{\tau}$ can be considered as upper bounds
(the actual effective reproduction numbers might be even lower).
Therefore, a feasible process determines in some sense the best use of given resources to achieve
the ``guaranteed lowest" number of infectious cases.

\section{NUMERICAL EXPERIMENTS}\label{example}

In this section we provide an example on a synthetic data set to demonstrate how the problems formulated above
can be used for controlling the spread of infection. In this example there are two regions ($m = 2$)
and for the sake of simplicity we assume that there is no transmission of infection between these regions
(i.e. $\beta_{1,2} = \beta_{2,1} = 0).$ Moreover, we consider only Problem 1; that is,
costs related to bed building ($c_{B}$), services ($c_{S}$) and before isolation ($c_{I}$) assumed to be sufficient in all cases.\\

{\bf Initial data.~} We assume the time interval is $[1,T] = [1,150];$ where $[1,100]$ is an initial (past) period and $[101,150]$ is the future/planning period for our optimal control problem.

The number of initial (i.e. $t \le 0$) infected cases is 2 in both regions. The set of possible values for
$\tau$ is $\{3,4,5\}$ (as in the case of data fitting in Section \ref{fit}).
We generate synthetic data - $x_1(t), x_2(t)$ for $t \in [1,100]$ by setting
\begin{itemize}
  \item Region 1: $\alpha = 0.6,$ $\beta_{11} = 0.30$ and $\tau_1(t) = 4, \forall t \in [1,50],$ $\tau_1(t) = 5, \forall t \in [51,100];$
  \item Region 2: $\alpha = 0.6,$ $\beta_{22} = 0.28$ and $\tau_2(t) = 4, \forall t \in [1,50],$ $\tau_2(t) = 5, \forall t \in [51,100].$
\end{itemize}

As in \cite{F1000}, $\omega$ as a gamma distribution function with mean value 7.5 defined by (\ref{gamma}). According to
formula (\ref{R}) corresponding effective reproduction numbers are
\begin{itemize}
  \item Region 1: $R = 0.90, 1.20$ and $1.48$ for $\tau = 3, 4$ and $5$, respectively;
  \item Region 2: $R = 0.84, 1.12$ and $1.38$ for $\tau = 3, 4$ and $5$, respectively.
\end{itemize}
Thus, the initial functions $x_1(t), x_2(t)$ for $t \in [1,100]$ have effective reproduction numbers
$R = 1.20$ and $1.48$ for $x_1(t)$ on $[1,50]$ and $[51,100]$, respectively;
$R = 1.12$ and $1.38$ for $x_2(t)$ on $[1,50]$ and $[51,100]$, respectively.\\

We will consider Problem 1 on the interval $[101,150].$
Both controls $\tau_i(t)$ will be assumed to be constant:
$\tau_i(t) = \tau_{i}, \forall t \in [101,150],$ $~ i =1,2.$
Values $\tau_1$ and $\tau_2$ will be used to describe future possible scenarios.

The initial number of beds are $b_1(100) = 126$ and $b_2(100) = 60.$
We assume that new beds will be created at times $T_1^b = 101,$  $T_2^b = 108,$ $T_3^b = 115$ and $T_4^b = 122;$
corresponding numbers of additional beds will be denoted by $\Delta b_i$, $i = 1, 2,3,4.$
Optimal control problem aims to distribute these additional beds between the regions.

We introduce a new variable - $\lambda_i$ that denotes part of $\Delta b_i$ considered for the first region, the remaining part $(1-\lambda_i) \Delta b_i$ for the second region.

Under these assumptions, Problem 1 can formulated as follows:
 $$
  {\rm Minimize_{(\lambda_1,\lambda_2,\lambda_3,\lambda_4)}~}~~  x_1(150) + x_2(150) $$
$$  + K \cdot \sum_{t = 101}^{150} \left[ \frac{h_1(t: \mathbf{\tau_1},\mathbf{x_1})}{b_1(t)} + \frac{h_2(t: \mathbf{\tau_2},\mathbf{x_2})}{b_2(t)} \right]
$$
$$
{\rm ~subject~to:~} \lambda_k \in [0,1], k = 1,2,3,4;
$$

$$
x_r(t+1) = \beta_{rr}  \, \sum_{i = 0}^{\tau_r(t)-1} (1-\alpha \omega(i)) \, x_r(t-d - i), r=1,2;
$$
$$
b_1(t) = b_1(100) + \sum_{i=1}^{|\{T_j^b \le t;~ j = 1,2,3,4\}|} \lambda_i \cdot \Delta b_i;
$$
$$
b_2(t) = b_2(100) + \sum_{i=1}^{|\{T_j^b \le t;~ j = 1,2,3,4\}|} (1-\lambda_i) \cdot \Delta b_i
$$
Here $K = 100$ and  $h_r(t: \mathbf{\tau_r},\mathbf{x_r})$ is calculated by (\ref{H}) by setting
$\sigma = 6$ and $d = 6.$\\

We can consider different scenarios depending on $\tau_1, \tau_2$ and additional beds $\Delta b_k,$ $k=1,2,3,4.$
Below we will provide three scenarios where we set
$$\Delta b_1 = 350, \Delta b_2 = 300, \Delta b_3 = 100, \Delta b_4 = 20$$ and change the values of $\tau_1, \tau_2$ as 3, 4 and 5 (that is, all possible values used in the data fitting problem in Section \ref{fit}). The aim here is to compare corresponding optimal distributions of bed capacities. We recall that additional beds are introduced (weekly) at
$T_1^b = 101,$  $T_2^b = 108,$ $T_3^b = 115$ and $T_4^b = 122.$
The cumulative number of infected cases at the start $t = 100$ are:  1259 in Region 1 and 675 in Region 2.\\

{\bf Case 1:~~} $\tau_1 = \tau_2 =3.$\\

The summary of the optimal solution obtained is provided below.
\begin{itemize}
  \item The optimal distribution of additional beds:
  $$
  \begin{array}{ccccc}
    Region 1: & 192.3 &  183.4 &  72.9 &  20 \\
    Region 2: & 157.7 & 116.6 & 27.1 &   0 \\ \hline
    Total:    &  350 & 300 & 100 & 20 \\
  \end{array}
$$
  \item Cumulative number of infected cases at the end $(t =150)$ of planning period: 2951 in Region 1 and 1289 in Region 2.
  \item The average number of bed occupancy over the time interval $[100,150]$ is
0.45 for Region 1 and 0.29 for Region 2.
  \item The maximum occupancy rates are: 0.83 (that is, average 0.83 patient per bed) in Region 1 and 0.55 in Region 2; that is,
the demand for hospital beds is met at every point $t\in [100,150].$
\end{itemize}

Thus the solution obtained is feasible.

\bigskip
{\bf Case 2:~~} $\tau_1 = \tau_2 = 4.$\\

Optimal solution obtained and some relevant parameters are:
\begin{itemize}
  \item The optimal distribution of additional beds:
  $$
  \begin{array}{ccccc}
    Region 1: & 192.0 & 183.4  & 74.5 &  20 \\
    Region 2: & 158.0 & 116.6 &  25.5 &   0 \\ \hline
    Total:    &  350 & 300 & 100 & 20 \\
  \end{array}
$$
  \item Cumulative number of infected cases at the end $(t =150)$ of planning period:
  5846 in Region 1 and 2259 in Region 2.
   \item The average number of bed occupancy over the time interval $[100,150]$ is
0.65 for Region 1 and 0.40 for Region 2.
    \item The maximum occupancy rates are: 0.96 (that is, average 0.96 patient per bed) in Region 1 and 0.52 in Region 2.
\end{itemize}

Again, the demand for hospital beds is met at every point $t\in [100,150]$ and accordingly the optimal solution obtained is feasible.\\

{\bf Case 3:~~} $\tau_1 = \tau_2 = 5.$\\

Optimal solution obtained and some relevant parameters are:
\begin{itemize}
  \item The optimal distribution of additional beds:
  $$
  \begin{array}{ccccc}
    Region 1: & 191.6 & 183.3 &  75 &  20 \\
    Region 2: & 158.4 & 116.7 &  25 &   0 \\ \hline
    Total:    &  350 & 300 & 100 & 20 \\
  \end{array}
$$
  \item Cumulative number of infected cases at the end $(t =150)$ of planning period:
  11585 (Region 1) and 4163 (Region 2).
  \item The average number of bed occupancy over the time interval $[100,150]$ is
0.83 for Region 1 and 0.51 for Region 2.
  \item The maximum occupancy rates are: 1.91 (that is, 1.91 patient per bed) in Region 1 and 1.05 in Region 2.
\end{itemize}

Therefore this solution is not feasible as the bed capacities are not enough for isolation all infected individuals.\\

Comparting the results in Cases 1-3, where $\tau_1 = \tau_2,$ we observe that, the optimal distributions of bed capacities are almost the same although the solution obtained in Case 3 is even not feasible.
In table \ref{tab2} we also provide the results obtained by assuming that the rate of increase in one region is greater than the other one. The results for $\tau_1 = \tau_2 + 1$ (that is, $\tau_1=3, \tau_2=4$ and $\tau_1=4, \tau_2=5$)
display quite similar optimal bed distributions.
We observe the same situation for $\tau_2 = \tau_1 + 1.$

Summarizing these results we note that the optimal control problem considered can provide quite ``robust" optimal distributions of new bed capacities across the regions under each of the assumptions $\tau_1 = \tau_2,$ $\tau_1 > \tau_2$ and $\tau_1 < \tau_2.$

\begin{table}
\caption{The optimal distribution of additional beds}\label{tab2}
\medskip
\centering
\begin{tabular}{|cc|ccccc|}
  \hline
 $\tau_1$ & $\tau_2$ &  &  &  &  &\\
\hline
 3 & 4 & Region 1: & 204.7 & 183.9 &  9.1 &  0 \\ 
   &   & Region 2: & 145.3 & 116.1 & 90.9 & 20 \\ \hline
 4 & 5 & Region 1: & 206.2 & 188.2 & 31.6 &  0 \\
   &   & Region 2: & 143.8 & 111.8 & 68.4 & 20 \\ \hline
 4 & 3 & Region 1: & 179.3 & 223.3 & 100 &  20 \\
   &   & Region 2: & 170.7 &  76.7 &   0 &   0 \\ \hline
 5 & 4 & Region 1: & 177.1 & 208.5 & 100 &  20 \\
   &   & Region 2: & 172.9 &  91.5 &   0 &   0 \\                    
 \hline
\end{tabular}
\end{table}

\end{document}